\documentclass[preprint,review,12pt]{elsarticle}

\usepackage{amssymb,amsmath,amsthm}

\newtheorem{dfn}{Definition}[section]
\newtheorem{thm}{Theorem}[section]

\journal{}

\begin{document}

\begin{frontmatter}



\title{A further generalization of random self-decomposability}

\author[label1]{S. Satheesh\corref{c1}}
\cortext[c1]{Corresponding author e-mail: \textit{ssatheesh1963@yahoo.co.in}}

\author[label2]{E. Sandhya\corref{c2}}
\cortext[c2]{e-mail: \textit{esandhya@hotmail.com}}

\address[label1]{NEELOLPALAM, S. N. Park Road, Trichur-680 004, India.}
\address[label2]{Department of Statistics, Prajyoti Niketan College, Pudukkad, Trichur-680\,301, India.} 

\begin{abstract}
The notion of random self-decomposability is generalized further. The notion is then extended to non-negative integer-valued distributions.
\end{abstract}


\end{frontmatter}

\section{Introduction}

Recently the notion of random self-decomposability ($RSD$) has been introduced by Kozubowski and Podg\'orski \cite{kp10} generalizing $SD$. They showed that if a CF is $RSD$ then it is both $SD$ and geometrically infinitely divisible ($GID$). Satheesh and Sandhya \cite{ss10} generalized this notion to Harris-RSD ($HRSD$) and showed that if a CF is $HRSD$ then it is both $SD$ and Harris-ID ($HID$). With this nomenclature $RSD$ is geometric-RSD ($GRSD$). Here we explore further generalizations of $HRSD$ \textit{viz.} $\mathcal{N}RSD$ and $\varphi RSD$, motivated by the elegent Proposition 2.3 in Kozubowski and Podg\'orski \cite{kp10}.
   
We need the notion of $\mathcal{N}$-infinitely divisible ($\mathcal{N}ID$) laws here. Let $\varphi$ be a Laplace transform (LT) that is also a standard solution to the Poincare equation, $\varphi(t)=P(\varphi(\theta t)), \theta \in \Theta$ where $P$ is a probability generating function (PGF) (\textit{see} Gnedenko and Korolev, \cite{gk96}, p.140).
\begin{dfn}{\rm
Let $\varphi$ be a standard solution to the Poincare equation and $N_\theta$ a positive integer-valued random variable (\textit {r.v.}) having finite mean with PGF $P_\theta (s)=
\varphi(\frac1\theta \varphi^{-1}(s))$, $\theta \in \Theta \subseteq (0,1)$. A characteristic function (CF) $f(t)$ is
$\mathcal{N}ID$ if for each $\theta \in \Theta$ there exists a CF $f_\theta(t)$ that is independent of $N_\theta$ such that $f(t) = P_\theta (f_\theta(t))$, for all $t \in \mathbf{R}$.
}\end{dfn}

\begin{thm}
(Gnedenko and Korolev, 1996, Theorem 4.6.3 on p.147) \cite{gk96} Let $\varphi$ be a standard solution to the Poincare equation. A CF $f(t)$ is $\mathcal{N}ID$ iff it admits the representation $f(t)=\varphi(-\log h(t))$ where $h(t)$ is a CF that is ID. $f(t)$ is $\mathcal{N}$stable if $h(t)$ is stable (p.151, \cite{gk96}).  
\end{thm} 

In the next section we describe $\mathcal{N}RSD$ laws and its discrete analogue in Section 3. In Section 4 we describe $\varphi RSD$ laws and its discrete analogue.

\section {$\mathcal{N}RSD$ distributions}

\begin{dfn}{\rm
A CF $f(t)$ is $\mathcal{N}RSD$ if for each $c \in (0,1]$ and each $\theta \in [0,1)$
\begin{equation}\label{e1.1}
f_{c,\theta}(t)=f_{c}(t). f_{\theta}(ct)
\end{equation}
is a CF, where $f_{c}(t)$ and $f_{\theta}(t)$ are given by
\begin {equation}
f_{c}(t)= \frac {f(t)}{f(ct)}
\end{equation}
\begin {equation}
f_{\theta}(t)=\varphi\{\theta\varphi^{-1}(f(t))\},   
\end{equation}
$\varphi$ being a standard solution to the Poincare equation. 
}\end{dfn} 

We now notice that the discussion leading to conceiving and proving Proposition 2.3 in Kozubowski and Podg\'orski \cite {kp10} holds in this generalization as well. When $c=1$ equation (1) becomes
\begin {equation}
f_{1,\theta}(t)=f_{\theta}(t)=\varphi\{\theta\varphi^{-1}(f(t))\}
\end{equation}
Or
\begin{equation}
f(t)=\varphi\{\frac1\theta\varphi^{-1}(f_{\theta}(t))\}
\end{equation}
for each $\theta \in [0,1)$. That is $f(t)$ is $\mathcal{N}ID$ and hence has no real zeroes. On the other hand since $\varphi(0)=1$, when $\theta=0$ equation (1) implies 
\begin {equation}
f_{c,0}(t)=f_{c}(t)=\frac{f(t)}{f(ct)} 
\end{equation}
is a CF for each $c \in (0,1]$. That is $f(t)$ is SD.

Conversely, if $f(t)$ is $SD$ then for each $c\in(0,1]$ the function $f_{c}(t)$ in (2) is a genuine CF and similarly if $f(t)$ is $\mathcal{N}ID$ then for each $\theta\in[0,1)$ the function $f_{\theta}(t)$ in (5) also is a genuine CF. Consequently (1) is a well defined CF.

\noindent \textbf{Remark 2.1} It may be noted that for the CF $f(t)$ to be $SD$ we only require that (a result due to Biggins and Shanbhag \textit see Fosum \cite{fo95}) (2) holds for all $c$ in some left neighbourhood of 1. Thus we may simplify the requirement here as: A CF $f(t)$ is $\mathcal{N}RSD$ if for each $c \in (a,1]$, and each $\theta \in [0,1)$ (1) holds, where $0<a<1$.

\noindent \textbf{Remark 2.2} In fact we may have apparently still weaker requirement in describing CFs that are $\mathcal{N}RSD$ as follows. A CF $f(t)$ is $\mathcal{N}RSD$ if for each $c \in (a,1)$, and each $\theta \in (0,1)$ (1) holds, where $0<a<1$. Now letting $c \uparrow 1$ we have $f(t)$ is $\mathcal{N}ID$. On the other hand letting $\theta \downarrow 0$ we have $f(t)$ is $SD$ since $lim_{\theta\downarrow 0}f_{\theta}(t)=1$, \textit{see e.g} Gnedenko and Korolev \cite{gk96}, page 149.

\noindent \textbf{Example 2.1} For the LT $\varphi(s)=(1+s)^{-\alpha}, \alpha>0$, $\varphi(\varphi^{-1}(s)/p)$ is a PGF of a non-degenerate distribution only if $\alpha=\frac{1}{k}, k\geq1$ integer, see Example 1 in Bunge \cite{bj96} or Corollary 4.5 in Satheesh \textit{et al.} \cite{sns02}. This PGF is that of Harris distribution (Satheesh et al. \cite {ssl10}) and the corresponding $\mathcal{N}RSD$ distribution is $HRSD$. When $k=1$ above, we have $GRSD$ ($RSD$ distributions of Kozubowski and Podg\'orski \cite {kp10}).
 
\noindent \textbf{Example 2.2} Invoking Theorem 1.1 when $\varphi(s)$ is $SD$ and $\log h(t)=-\lambda|t|^\alpha$ we have, for each $c\in(a,1]$ 
\begin {equation}
f(t)=\varphi(|t|^\alpha)=\varphi(c|t|^\alpha).\varphi_c(|t|^\alpha).
\end{equation}
That is $f(t)$ is both $SD$ and $\mathcal{N}$-strictly stable. Thus we have a good collection of CFs that are both $SD$ and $HID$ and thus $HRSD$. Kozubowski and Podg\'orski \cite {kp10}) present examples of a variety of CFs $h(t)$ that are stable.  

\section {Discrete analogue of $\mathcal{N}RSD$ distributions}
Steutel and van Harn \cite{sh79} had described discrete SD ($DSD$) distributions. Satheesh and Sandhya \cite{ss10} have described $DHRSD$, discrete analogue of $HRSD$ distributions. We now introduce discrete $\mathcal{N}RSD$ ($D\mathcal{N}RSD$) distributions.

\begin{dfn}{\rm
(Satheesh \textit{et al.} \cite{ssl10}) Let $\varphi$ be a standard solution to the Poincare equation and $N_\theta$ a positive integer-valued \textit {r.v.} having finite mean with PGF $P_\theta (s)=
\varphi(\frac1\theta \varphi^{-1}(s))$, $\theta \in \Theta \subseteq (0,1)$. A PGF $P(s)$ is
$D\mathcal{N}ID$ if for each $\theta \in \Theta$ there exists a PGF $Q_\theta (s)$ that is independent of $N_\theta$ such that $P(s) = P_\theta (Q_\theta(s))$, for all $|s|\leq1$.
}\end{dfn}

\begin{thm}
(Satheesh \textit{et al.} \cite{ssl10}) Let $\varphi$ be a standard solution to the Poincare equation. A PGF $P(s)$ is $D\mathcal{N}ID$ iff it admits the representation $P(s)=\varphi(-\log R(s))$ where $R(s)$ is a PGF that is DID.
\end{thm} 

\begin{dfn}{\rm
A PGF $P(s)$ is $D\mathcal{N}RSD$ if for each $c \in (0,1]$ and each $\theta \in [0,1)$
\begin{equation}\label{e1.1}
P_{c,\theta}(s)=P_{c}(s). Q_{\theta}(1-c+cs)
\end{equation}
is a PGF, where $P_{c}(s)$ and $Q_{\theta}(s)$ are given by
\begin {equation}
P_{c}(s)= \frac {P(s)}{P(1-c+cs)}
\end{equation}
\begin {equation}
Q_{\theta}(s)=\varphi\{\theta\varphi^{-1}(P(s))\},   
\end{equation}
$\varphi$ being a standard solution to the Poincare equation. 
}\end{dfn}

We may now proceed as in Section 2 describing the relation between $DSD$, $D\mathcal{N}ID$ and $D\mathcal{N}RSD$ distributions. Further, remarks similar to Remarks 2.1 and 2.2 are relevant here also and Examples on the lines of Example 2.1 nad 2.2 can also be discussed.

\section {$\varphi RSD$ distributions}
A further generalization of $\mathcal{N}RSD$ distributions is possible invoking the notion of $\varphi ID$ law that generalizes $\mathcal{N}ID$ laws, see Satheesh \cite{ss04} and Satheesh \textit{et al.} \cite{ssl10}, \cite{ssl11}) for its discrete analogue. We first describe the discrete case.

\begin{dfn}{\rm
(Satheesh \textit{et al.} \cite{ssl10}) Let $\varphi$ be a LT. A PGF $P(s)$ is $D\varphi ID$ if there exists a sequence $\{\theta_n\} \downarrow 0$ as $n\rightarrow\infty$ and a sequence of $PGF$s $Q_n(s)$ such that
\begin {equation}
P(s)=lim_{n\rightarrow\infty}\varphi(\frac{1-Q_n(s)}{\theta_n}).
\end{equation}
}\end{dfn}

\begin{thm}
(Satheesh \textit{et al.} \cite{ssl11}) Let $\{Q_{\theta}(s), \theta \in \Theta\}$ be a family of PGFs and $\varphi$ a LT. Then 
\begin {equation}
lim_{\theta \downarrow 0}\varphi(\frac{1-Q_{\theta}(s)}{\theta})
\end{equation}
exists and is $D\varphi ID$ iff there exists a PGF $R(s)$ that is DID such that
\begin {equation}
lim_{\theta \downarrow 0}\frac{1-Q_{\theta}(s)}{\theta}=-\log R(s)
\end{equation} 
\end{thm}

\begin{dfn}{\rm
A PGF $P(s)$ is $D\varphi RSD$ if for each $c \in (a,1)$ and each $\theta \in (0,b), 0<a,b<1$
\begin{equation}\label{e1.1}
P_{c,\theta}(s)=P_{c}(s). Q_{\theta}(1-c+cs)
\end{equation}
is a PGF, where $P_{c}(s)$ and $Q_{\theta}(s)$ are given by
\begin {equation}
P_{c}(s)= \frac {P(s)}{P(1-c+cs)}
\end{equation}
\begin {equation}
Q_{\theta}(s)=1- \theta \varphi^{-1}(P(s)).   
\end{equation}
}\end{dfn}

The restriction of $\alpha=\frac{1}{k}, k\geq1$ integer in Example 2.1 is not in this notion. We may now proceed as in Section 3 describing the relation between $DSD$, $D\varphi ID$ and $D\varphi RSD$ distributions. This has been possible since $lim_{\theta\downarrow 0}Q_{\theta}(t)=1$. The case of $\varphi RSD$ follows on similar lines.


\end{document}